\documentclass[11pt]{article}

\usepackage{amssymb,amsmath,amsthm,rotating}

\newcommand{\bb}[1]{\mathbb{#1}}
\newcommand{\cl}[1]{\mathcal{#1}}

\newcommand{\sst}{\scriptstyle}

\theoremstyle{plain}
\newtheorem{lem}{Lemma}[section]

\newtheorem{tm}{Theorem}

\theoremstyle{definition}

\begin{document}

\title{Essentially Reductive Hilbert Modules II}

\author{Ronald G.~Douglas}

\date{}
\maketitle

\abstract{
Many Hilbert modules over the polynomial ring in $m$ variables are essentially reductive, that is, have commutators which are compact. Arveson has raised the question of whether the closure of homogeneous ideals inherit this property and provided motivation to seek an affirmative answer. Positive results have been obtained by Arveson, Guo, Wang and the author. More recently, Guo and Wang extended the result to quasi-homogeneous ideals in two variables. Building on their techniques, in this note the author extends this result to Hilbert modules over certain Reinhardt domains such as ellipsoids in two variables and analyzes extending the result to the closure of quasi-homogeneous ideals in $m$ variables when the zero variety has dimension one.}

\setcounter{section}{-1}

\section{Introduction}\label{sec0}

\indent

In his study \cite{Ar3} of the $m$-shift Hilbert module $H^2_m$ over ${\bb C}[\pmb{z}]$, the polynomial algebra in $m$ variables, Arveson formulated a provocative conjecture  which has attracted the attention of several researchers. He had established that the commutators of all operators and their adjoints defined on $H^2_m$ by module multiplication by polynomials in ${\bb C}[\pmb{z}]$ belonged to the Schatten $p$-class ${\cl L}^p$ for all $p>m$. He then conjectured that the same commutator condition held for all submodules in $H^2_m$ obtained as the closure $[I]$ in $H^2_m$ of a homogeneous ideal $I$ and established the result in \cite{Ar4} in case $I$ is generated by monomials.\footnote{The same result was proved earlier for the quotient defined by every homogeneous submodule in the Hardy space $H^2({\bb D}^2)$ for the bidisk ${\bb D}^2$ by Curto, Muhly and Yan \cite{CMY}.}$^{,}$\footnote{Such Hilbert modules were defined to be essentially reductive in \cite{DP} and, later, essentially normal in \cite{Ar4}.\par This research was prompted by discussions with K.~Guo and K.~Wang during a visit to Fudan University in July, 2005 and with Arveson and other researchers during a DST-NSF funded visits to Bangalore in December, 2003 and 2005.\par
Mathematics Subject Classification:\ 32T15, 46L80, 46M20, 47L15.\par
Key words and phrases:\ Hilbert modules, submodules, essentially reductive, essentially normal, $K$-homology.} Using somewhat different methods, the author extended this latter result in \cite{Dou2} to a larger family of commuting weighted shifts in $m$-variables. At the same time, it was pointed out that these results suggested extending the conjecture to domains other than the unit ball ${\bb B}^m$ in ${\bb C}^m$. The meaning of this latter statement becomes clear once one understands that $H^2_m$ can be realized as a reproducing Hilbert space of holomorphic functions on ${\bb B}^m$ and that the commutator condition on $H^2_m$ can be shown to be equivalent to the same statement for the Bergman space $L^2_a({\bb B}^m)$ on ${\bb B}^m$, that is, for the closure of ${\bb C}[\pmb{z}]$ in the $L^2$-space on ${\bb B}^m$ relative to volume measure.

Subsequently, Guo obtained results in \cite{G} which, when combined with earlier techniques of Arveson, established Arveson's original conjecture in $H^2_2$ for arbitrary homogeneous ideals in ${\bb C}[z_1,z_2]$. Then Guo, joined by Wang in \cite{GW}, extended this result  to general $m$ when the complex dimension of the zero variety $Z(I)$ of $I$ is one or, equivalently, when the Hilbert polynomial for the quotient module $H^2_m/[I]$ is linear (cf.\ \cite{DY}). Finally, in a  recent paper \cite{GW2} Guo and Wang showed that the closure of quasi-homogeneous ideals in a class of Hilbert modules, including $H^2_2$, are essentially reductive.

In this note, based largely on the methods of Guo and Wang, we extend this latter result 
in two ways. First, we show that submodules obtained as the closure of quasi-homogeneous ideals in the Bergman space for certain Reinhardt domains in ${\bb C}^2$ are essentially reductive. Note that this result includes the case of the closure of homogeneous ideals, which is new. Second, we show that the same is true for the closure of certain, very restrictive quasi-homogeneous ideals in ${\bb C}[\pmb{z}]$ for arbitrary $m$ with one dimensional zero variety in the Bergman space for nice Reinhardt domains in ${\bb C}^m$. This collection of Reinhardt domains  includes the ellipsoids of the form $E_{\pmb{a}} = \left\{\pmb{z}\in {\bb C}^m\colon\ \sum\limits^m_{i=1} a_i|z_i|^2 < 1\right\}$ for $\pmb{a}$ in $[0,\infty)^m$. Of course, $E_{\pmb{1}} = {\bb B}^m$ with $\pmb{1} = (1,\ldots, 1)$.

We will say that an ideal in ${\bb C}^m$ is bivariate if it is generated by polynomials in two of the variables $z_1,\ldots z_m$ at a time. Note that all ideals in ${\bb C}[z_1,z_2]$ are bivariate.  The limit of the techniques in this note would seem to be bivariate ideals with one dimensional zero variety although we are unable to establish such a result at this time.

Most of the results in \cite{Ar4}, \cite{Dou2}, \cite{G} and \cite{GW} apply not just to the closures of homogeneous ideals in a Hilbert space completion ${\cl H}$ of ${\bb C}[\pmb{z}]$ but to homogeneous submodules of ${\cl H}\otimes {\bb C}^k$. In this note, we confine our attention to the multiplicity one case, $k=1$, or the closure of ideals.

We will assume the reader is familiar with \cite{GW2} although we will provide statements of the relevant definitions, lemmas and propositions and the necessary proofs but emphasize mainly the  points that are different and not straightforward extensions of those in \cite{GW2}.

\section{The Basic Setup}\label{sec1}

\indent 

A Reinhardt domain $\Omega$ in ${\bb C}^m$ is one for which $\pmb{z} = (z_1,\ldots, z_m)$ in $\Omega$ implies $(e^{i\theta_1}z_1,\ldots, e^{i\theta_m}z_m)$ is in $\Omega$ for all $m$-tuples $(e^{i\theta_1},\ldots, e^{i\theta_m})$ in the $m$-torus. Hence, the absolute values of the coordinates of a point are sufficient to determine if the point is in $\Omega$. Even so, Reinhardt domains can be somewhat pathological. A nice class of them is determined by smooth functions $\varphi\colon \ [0,\infty)^m\to [0,\infty)$ so that $\Omega_\varphi = \{\pmb{z}\in {\bb C}^m\colon \ \varphi(|z_1|,\ldots, |z_m|) < 1\}$. We assume further that $\varphi$ is monotonically increasing in each variable $t_i$, $i=1\ldots m$. In this case, the boundary $\partial\Omega_\varphi$ of $\Omega_\varphi$ is the set $\{\pmb{z}\in {\bb C}^m\colon \ \varphi(|z_1|,\ldots, |z_m|) = 1$. Moreover, $\Omega_\varphi$ is homeomorphic to ${\bb B}^m$ and $\partial\Omega_\varphi$ to the unit sphere $\partial{\bb B}^m$.
 However, unlike the case of analogous domains in ${\bb C}$, there is, in general, no biholomorphism between one of  them and the ball. 

Now ${\bb C}[\pmb{z}]$ is norm dense in the closure, ${\cl H}_\varphi=L^2_a(\Omega_\varphi)$, of the functions holomorphic on a neighborhood of the closure of $\Omega_\varphi$ for volume measure on $\Omega_\varphi$. Moreover, since  volume measure on $\Omega_\varphi$ restricts, in the sense of the Fubini Theorem, to Lebesgue measure on the natural $m$-torus through each point, we see that the monomials form an orthogonal basis for ${\cl H}_\varphi$. Further,  $\Omega_\varphi$ is  pseudo-convex  if and only if it is logarithmically convex  (cf.\ \cite{R}). In addition, if $\partial\Omega_\varphi$ contains no disks (cf.\ \cite{FS}), then the commutators on ${\cl H}_\varphi$ are in ${\cl L}^p$ for $p>m$. (We will say that such Reinhardt domains are pseudo-convex without disks.) This implies that ${\cl H}_\varphi$ satisfies the hypotheses of Theorem 4.3 in \cite{Dou2}. Hence, modules defined as the closure of ideals in ${\bb C}[\pmb{z}]$ generated by monomials are $p$-essentially reductive for $p>m$. In this note, we show the same is true when the ideal is bivariate,  quasi-homogeneous, has zero variety with complex dimension one and  is \emph{radical}. Note that in case if the ideal is actually homogeneous, then for some Reinhardt domains this result follows from the earlier work of Guo and Wang \cite{GW}.

We assume in what follows that $\Omega_\varphi$ is a  pseudo-convex Reinhardt domain without disks. A polynomial $p(\pmb{z}) = \sum\limits^\infty_{\pmb{i},\pmb{j}=0} p_{\pmb{i}}\pmb{z}^{\pmb{i}}$, where $\pmb{i} = (i_1,\ldots, i_m)$ is said to be quasi-homogeneous of degree $\ell$  with respect to the weights $\pmb{n} = (n_1,\ldots, n_m)$ if $p_{\pmb{i}}\ne 0$ implies $\sum\limits^m_{i=1} n_ji_j=\ell$.

If we consider the action of ${\bb T}^m$ on $\Omega_\varphi$ in the obvious way, then quasi-homogeneity can be characterized in terms of this action. For each $\lambda$ in ${\bb R}$, consider the action $\gamma_\lambda$ on $\Omega_\varphi$ defined by
\[
\gamma_\lambda(z_1,\ldots, z_m) = (e^{in_1\lambda}z_1,\ldots, e^{in_m\lambda}z_m).
\]
Then the quasi-homogeneous polynomials of degree $\ell$ for the weights $(\pmb{n})$ are just the eigenvectors for the action of $\{\gamma_\lambda\}_{\lambda\in{\bb R}}$ on ${\cl H}_\varphi$ for the eigenvalue $\ell$. Since the monomials form an orthogonal basis for ${\cl H}_\varphi$, we see for fixed weights $(\pmb{n})$ that ${\cl H}_\varphi = \oplus {\cl H}^\ell_\varphi$, where the ${\cl H}^\ell_\varphi$ are the quasi-homogeneous polynomials of degree $\ell$ with respect to $(\pmb{n})$, and the decomposition is orthogonal. As indicated by Guo and Wang, if ${\cl M}$ is the closure of a quasi-homogeneous ideal with weights $(\pmb{n})$, then
\[
{\cl M} = \oplus {\cl M}_\ell \text{ and } {\cl M}^\bot = \oplus {\cl H}^\ell_\varphi \ominus {\cl M}_\ell,
\]
where ${\cl M}_\ell = {\cl M}\cap {\cl H}^\ell_\varphi$ and ${\cl H}^\ell_\varphi \ominus {\cl M}_\ell = {\cl M}^\bot \cap {\cl H}^\ell_\varphi$. In other words, both ${\cl M}$ and ${\cl M}^\bot$ are graded.

Let ${\cl M}$ be the closure in ${\cl H}_\varphi$ of a quasi-homogeneous ideal $I$. For $p$ in ${\bb C}[\pmb{z}]$, let $M_p$ denote the operator on ${\cl H}_\varphi$ defined by module multiplication by $p$; $A_p$ the restriction of $M_p$ to the submodule ${\cl M}$, $B_\varphi$ the action of $M_p$ from ${\cl M}$ to ${\cl M}^\bot$, and $C_p$ the compression of $M_p$ to ${\cl M}^\bot$, that is, $M_p = \left(\begin{smallmatrix} A_p&B_p\\ 0&C_p\end{smallmatrix}\right)$.

Using the index formula for Toeplitz operators on ${\cl H}_\varphi$ \cite{CMY}, one can show that there exists a Fredholm operator on ${\cl H}_\varphi \otimes {\bb C}^m$  with non-zero index. Hence, there exists a non-zero compact operator in the $C^*$-algebra, $C^*({\cl H}_\varphi)$, generated by the operators defined by module multiplication on ${\cl H}_\varphi$. Further, since ${\cl H}_\varphi$ is a reproducing kernel Hilbert space over $\Omega_\varphi$ which is connected, ${\cl H}_\varphi$ is irreducible. Therefore, the ideal ${\cl K}({\cl H}_\varphi)$ of all compact operators on ${\cl H}_\varphi$ is contained in $C^*({\cl H}_\varphi)$. Finally, the same index result implies that $C^*({\cl H}_\varphi)/{\cl K}({\cl H}_\varphi)$ is isomorphic to $C(\partial\Omega_\varphi)$ with the homomorphism extending the obvious map $M_p\to p|_{\partial \Omega_\varphi}$ for $p$ in ${\bb C}[\pmb{z}]$. (This is a known result but perhaps not quite in this generality \cite{C}.) Thus we have the short exact sequence
\[
0 \longrightarrow {\cl K}({\cl H}_\varphi) \to {\cl C}^*({\cl H}_\varphi) \overset{\sst\pi}{\longrightarrow} C(\partial\Omega_\varphi) \longrightarrow 0.
\]
where $\pi$ is the quotient map.

\section{Main Results}\label{sec2}

\indent

We recall a lemma from \cite{Ar4} and  \cite{Dou2} which is one starting point for Guo and Wang \cite{GW2}. For ${\cl M}$ a subspace of the Hilbert space ${\cl H}$, let $P_{\cl M}$ denote the projection onto ${\cl M}$.

\begin{lem}\label{lem2.1a}
If ${\cl H}$ is an essentially reductive Hilbert module over ${\bb C}[\pmb{z}]$ and ${\cl M}$ is a submodule of ${\cl H}$, then ${\cl M}$ is essentially reductive iff ${\cl H}/{\cl M}$ is essentially reductive iff the operators $P_{{\cl M}}M_pP_{{\cl M}^\bot}$ are compact for all $p$ in ${\bb C}[\pmb{z}]$.
\end{lem}

Next we recall another result from \cite{Dou2} which we use to establish an extension of Lemma 2.5 of \cite{GW2}.
The result in \cite{Dou2} is the dual of this one obtained by taking adjoints.

\begin{lem}\label{lem2.2a}
If ${\cl M}_0, {\cl M}_1$ and ${\cl M}_2$ are Hilbert modules over ${\bb C}[\pmb{z}]$ with ${\cl M}_1$ and ${\cl M}_2$ essentially reductive and $X_0$ and $X_1$ are module maps so that $X_0$ is isometric and
\[
0 \longrightarrow {\cl M}_0 \overset{\sst X_0}{\longrightarrow} {\cl M}_1 \overset{\sst X_1}{\longrightarrow} {\cl M}_2 \longrightarrow 0
\]
is exact, then ${\cl M}_0$ is essentially reductive.
\end{lem}

\begin{lem}\label{lem2.3a}
If $I$ and $J$ are ideals in ${\bb C}[\pmb{z}]$ so that $Z(I)\cap Z(J) \cap \partial \Omega_\varphi = \emptyset$ and both $[I]$ and $[J]$ are essentially reductive, then  $[I\cap J]$ is essentially reductive.
\end{lem}

\begin{proof}
We begin by first defining the isometric module map $X_0\colon \ [I\cap J] \longrightarrow [I]\oplus [J]$ so that $X_0f = \frac1{\sqrt 2}(f\oplus -f)$. Next we define $X_1\colon \ [I]\oplus [J] \longrightarrow [I,J]$ so that $X_1(f\oplus g) = f+g$ and observe that $0\longrightarrow [I\cap J] \overset{\sst X_0}{\longrightarrow} \mathop{\oplus}\limits^{[I]}_{[J]} \longrightarrow [I,J]\longrightarrow 0$ is seen to be exact once we know that $X_1$ is onto. To show that, let $\{p_i\}^s_{i=1}$ and $\{q_j\}^t_{j=1}$ be finite sets of generators for $I$ and $J$, respectively. If we consider the row operators $R$ and $S$ from ${\cl H}\oplus\cdots\oplus{\cl H}$ and ${\cl H}\oplus\cdots{\cl H}$ defined by
\[
R(f_1\oplus\cdots\oplus f_s) = \sum^s_{i=1} p_if_i \quad \text{and}\quad S(g_1\oplus\cdots g_t) = \sum^t_{j=1} q_jg_j,
\]
respectively, then $RR^* + SS^* = \sum\limits^s_{i=1} M_{p_i}M^*_{p_i} + \sum\limits^t_{j=1} M_{q_j}M^*_{q_j}$. Moreover, since $Z(I)\cap Z(J)\cap \partial\Omega_\varphi = \emptyset$, we have $\pi(RR^*+SS^*) = \sum\limits^s_{i=1} |p_i|^2 + \sum\limits^t_{j=1} |q_j|^2>0$ on $\partial\Omega_\varphi$ implying that $RR^*+SS^*$ is Fredholm. But the range of $X_1$ contains the span of the ranges of $R$ and $S$, which has finite codimension in ${\cl H}_\varphi$. Thus, $X_1$ is onto the essentially reductive model $[I,J]$.
\end{proof}

Note that one can modify the proof in \cite{GW2} to establish that $[I\cdot J]$ is essentially reductive under the assumptions of this lemma but we won't need that result in this paper.

We now follow \cite{GW2} in reducing the case of quasi-homogeneous submodules to principal ideals generated by  a power of a prime. We begin with the case of two variables.

If ${\cl M}$ is a quasi-homogeneous submodule of ${\cl H}_\varphi$, then $I = {\cl M}\cap {\bb C}[z_1,z_2]$ is a quasi-homogeneous ideal with $[I]={\cl M}$. (The argument to establish this result depends on decomposing elements of ${\cl M}$ into components by degree and then establishing that the components are also in ${\cl M}$.) The remainder of the proof rests on a result of Yang \cite{Y} which states that $I$ has a Beurling form, $I=pL$, where $p$ is the GCD of the polynomials in $I$ and $L$ is an ideal  with finite codimension in ${\bb C}[z,w]$. 
Since $L$ is quasi-homogeneous and ${\bb C}[z,w]/L$ is finite dimensional, one knows that $Z(L) = (0)$. Thus one can use Lemma \ref{lem2.3a} to conclude that $[p]\cap [L] = [pL]$ is essentially reductive if $[p]$ is.

 Since a quasi-homogeneous polynomial $p(z_1,z_2)$ has the form
$p(z_1,z_2) = \alpha z^r_1z^s_2\prod_i(z^{t_i}_1-\alpha_iz^{u_i}_2)^{v_i}$ with the $\alpha_i$ non-zero and distinct \cite{BM}, a repeated application of Lemma \ref{lem2.3a} reduces the essential reductivity of ${\cl M}$ to submodules generated by ideals of the form $[z^r_1z^s_2]$ and $[(z^t_1-\alpha z^u_2)^v]$.  For the first case, the result follows from the theorem for monomials in \cite{Dou2}. For the second, we need a further argument.

First, we extend Lemma 2.2 of Guo and Wang \cite{GW2} to a Bergman space ${\cl H}_\varphi$ over a  pseudo-convex Reinhardt domain $\Omega_\varphi$ in ${\bb C}^m$ without disks in its boundary. We divide the proof into two parts. For ${\cl H}$ a Hilbert  space, let $\pi$ denote the quotient map from ${\cl L}({\cl H})$ onto the Calkin algebra ${\cl L}({\cl H})/{\cl K}({\cl H})$. (This notation is consistent with the earlier definition of $\pi$.)

\begin{lem}\label{lem2.1}
If ${\cl A}_I$ is the $C^*$-subalgebra of $C^*({\cl H}_\varphi)$ generated by the operators $M_p$ for $p$ in an ideal $I$ in ${\bb C}[\pmb{z}]$ and ${\cl K}({\cl H}_\varphi)$, then
\[
{\cl A}_I = \{T\in C^*({\cl H}_\varphi)\colon  \pi(T) \equiv 0 \text{ on } Z(I)\}.
\]
\end{lem}

\begin{proof}
Since $\pi({\cl A}_I)$ is the self-adjoint subalgebra of $C(\partial \Omega_\varphi)$ generated by the restriction of the functions in $I$ to $\partial \Omega_\varphi$, the result follows.
\end{proof}

\begin{lem}\label{lem2.2}
Let $I$ be an ideal in ${\bb C}[\pmb{z}]$ so there exist positive constants $\pmb{a}=(a_1,\ldots, a_m)$ such that $M_{z_i}M^*_{z_i}-a_i$ is in ${\cl A}_I$ for $i=1,2,\ldots, m$.
If ${\cl M}$ is the closure of  $I$ in ${\bb C}[\pmb{z}]$ and $Z(I)\cap \partial \Omega_\varphi\cap Z(z_1) = \emptyset$, then $C_{z_1}C^*_{z_1}$ is Fredholm and $C_{z_1}B^*_{z_1}$ is compact.
\end{lem}

\begin{proof}
For $f$ in $I$, we see that the matrix  $M_f = \left(\begin{smallmatrix} A_f&B_f\\ 0&0\end{smallmatrix}\right)$, since the range of $M_f$ is contained in ${\cl M}$. But $M_f$ is essentially normal from which it follows that $B_f$ is compact for $f$ in $I$.

By the Hilbert basis theorem, there exist polynomials $q_1,\ldots, q_s$ in $I$ that generate $I$. Hence, we have $Z(I) = \{\pmb{z}\in {\bb C}^m\colon \ q_1(\pmb{z}) =\ldots= q_s(\pmb{z}) = 0\}$. Therefore, $|\pmb{z}_1|^2 + \sum\limits^s_{i=1} |q_i(\pmb{z})|^2>0$ on $\partial \Omega_\varphi$ and hence $T = M_{z_1}M^*_{z_1} + \sum\limits^s_{i=1} M_{q_i}M^*_{q_i}$ is Fredholm.

Now $[T,P_{\cl M}]$ compact implies that $C_{z_1}C^*_{z_1}$ is Fredholm and $C_{z_1}B^*_{z_1}$ is compact using matrix calculations as in \cite{GW2}.
\end{proof}

\begin{lem}\label{lem2.4}
If, in addition to the hypotheses of Lemma \ref{lem2.2}, either $m=2$ and $I$ is the principal ideal generated by $(z^t_1-\alpha z^u_2)^v$ with $\alpha\ne 0$ or $m$ is arbitrary and $I$ is radical, then $B_{z_1}$ is compact.
\end{lem}

\begin{proof}
 Since $C_{z_1}C^*_{z_1}$ is Fredholm, this will follow once one knows that $C_{z_1}$ is Fredholm or that the null space of $C_z$ is finite dimensional. We claim, in fact, that it's (0) in both cases.

If $C_{z_1}h=0$ for some $h$ in ${\cl M}^\bot$, then $z_1h$ is in ${\cl M}$. We can decompose $h$ in ${\cl H}_\varphi$ so that $h=\Sigma h_\ell$ with $h_\ell$ in ${\cl H}^\ell_\varphi$. But then $z_1h=\Sigma z_1h_\ell$ is in ${\cl M}$ and since $I$ is quasi-homogeneous, we see that $z_1h_\ell$ is in ${\cl M}_{\ell+1} \subset I\cap {\cl M} = I$. Hence, we have $Z(z_1)\cup Z(h_\ell) =  Z(z_1h_\ell) \supset  Z(I)$ and since $Z(z_1) \cap  Z(I)= \emptyset$, it follows that $Z(I)\subset  Z(h_\ell)$. Since $I$ is radical, we have $h_\ell$ in ${\cl M}$, which implies $h_\ell=0$ or $h=0$ and the result is proved in this case.

If $I$ is singly generated by $q(z_1,z_2)=(z^t_1-\alpha z^u_2)^v$, then  we obtain $z_1h_\ell = qp$. Since $z_1$ and $q$ are relatively prime, it follows that $z_1$ divides $p$ and hence $h_\ell$ is in ${\cl M}_\ell$ and is 0. This is the argument in \cite{GW2}. 
\end{proof}

We complete the proof of the main theorem in essentially the same way as in \cite{GW2}.

\begin{tm}\label{thm1}
If $\Omega_\varphi$ is a  pseudo-convex Reinhardt domain in ${\bb C}^2$ without disks in its boundary, then every quasi-homogeneous submodule of ${\cl H}_\varphi$ is essentially reductive.
\end{tm}

\begin{proof}
The first requirement needed to apply Lemma \ref{lem2.4} is to show somehow that $|z_1|$ is constant on $\partial\Omega_\varphi\cap Z(I)$. Since $I$ is generated by $(z^t_1-\alpha z^u_2)^v$, we see that there exists a monotonically increasing function $\psi\colon \ [0,\infty)\to [0,\infty)$ so that $(z_1,z_2)$ in $Z(I)$ implies $\psi(|z_1|) = |z_2|$. Thus if $(z_1,z_2)$ is in $\partial\Omega_\varphi\cap Z(I)$ we have $\varphi(|z_1|, \psi(|z_1|))=1$. But this uniquely determines $|z_1|$ since $\varphi$ is monotone in each variable. We can now apply Lemma \ref{lem2.4} to conclude that $B_{z_1}$ is compact. Reversing the roles of $z_1$ and $z_2$ yields that $B_{z_2}$ is compact which completes the proof.
\end{proof}

What kind of assumption can we make for $m>2$ to establish the hypothesis of the constancy of the restriction of $|z_1|$ to $\partial\Omega_\varphi\cap Z(I)$? One possibility is the following definition.

An ideal $I$ in ${\bb C}[\pmb{z}]$ will be  said to have an absolutely determining zero variety if for some fixed $j_0$, $1\le j_0\le m$, and each $i$, $1\le i\le m$, there is a continuous monotonically increasing function $\psi_i$ so that $|z_i| = \psi_i(|z_{j_0}|)$ for each $\pmb{z}$ in $Z(I)$. 

Note that one can show that if $I$ is an ideal in ${\bb C}[z_1,z_2]$  generated by a quasi-homogeneous polynomial of the form $p(z_1,z_2)= (z^t_1-\alpha z^u_2)^v$ for nonnegative integers $t,u$ and $v$, and $\alpha$ in ${\bb C}\backslash\{\pmb{0}\}$, then $I$ has an absolutely determining zero variety. Also, if $I$ is an ideal in ${\bb C}[z_1,\ldots, z_m]$  generated by such a polynomial $q$ in $z_1$ and $z_2$ and the monomials $z_3,z_4,\ldots, z_m$, then the same thing is true. Further, if $q_1(z_1,z_m)$ and $q_2(z_2,z_m),\ldots, q_k(z_k,z_m)$ for $1\le k<m-1$ are each quasi-homogeneous polynomials of the same form and $I$ is the ideal generated by them and the monomials $z_{k+1},\ldots, z_{m-1}$, then again $I$ has an absolutely determining zero variety. Finally, other combinations are possible such as an ideal generated by quasi-homogeneous polynomials of the above form $p_1(z_1,z_2)$, $p_2(z_2,z_3)$ and $p_3(z_3,z_4)$ in ${\bb C}[z_1,\ldots, z_4]$. 

Continuing this line of thought one can show  for $I$ a bivariate, quasi-homogeneous, radical ideal having $\dim Z(I)=1$, that $I$ has an absolutely determining zero variety. Note that the fact that $I$ is radical forces the generating polynomials to  be prime factors having the form $z^t_i -\alpha z^u_j$ for $\alpha\ne 0$. The fact that $\dim Z(I)=1$ forces the absolute values to all be related to that of one of the coordinates  $z_i$. Let us provide some more details.

Let $I$ be a bivariate ideal and set $I_i = I\cap {\bb C}[z_i]$ for $1\le i\le m$ and $I_{i,j} = I\cap {\bb C}[z_i,z_j]$ for $1\le i,j\le m$, $i\ne j$. (Here, we are viewing ${\bb C}[z_i]$ and ${\bb C}[z_i,z_j]$ as subalgebras of ${\bb C}[\pmb{z}]$ in the obvious way.) Since $I_i\subset I_{i,j}$ for all $1\le j\le m$, $j\ne i$, $I$ is generated by the collection $[I_{i,j}]$, and we have $Z(I)=\cap Z(I_{i,j})$. Moreover, if $I$ is a radical ideal, then so is each $I_i$ and $I_{i,j}$. Finally, if $I$ is quasi-homogeneous for the weights $(\pmb{n})$, then $I_{i,j}$ is quasi-homogeneous for the weights $(n_i,n_j)$.

If $I_i\ne (0)$, then it is generated by $z_i$ and hence $z_i$ is in $I$. If $z_i$ is in $I$ and $I_{i,j}\ne I_i$, then $I_j\ne (0)$ and $I_{i,j}$ is generated by $z_i$ and $z_j$. Moreover, both $z_i$ and $z_j$ are in $I$. Partition the integers $\{1,2,\ldots, m\}$ into two sets $\Gamma_1$ and $\Gamma_2$ so that $i$ is in $\Gamma_1$ if and only if $z_i$ is in $I$. Then ${\bb C}[z_i\colon \ i\in\Gamma_1]\subset I$ and $Z(I)\subset \{\pmb{z}\in {\bb C}^m\colon \ z_i=0$ for $i\in\Gamma_1\}$. We seek now to partition $\Gamma_2$ so that $i$ and $j$ in $\Gamma_2$ are equivalent if $I_{i,j}\ne (0)$. Since the ideal $I_{i,j}\subset {\bb C}[z_i,z_j]$ is generated by a prime polynomial $z^s_i-\alpha z^t_j$ for positive integers $s$ and $t$ and $\alpha\ne 0$, we see that we obtain an equivalence relation which partitions $\Gamma_2$ into subsets $\{\Gamma^k_2\}$ of $\{1,\ldots, m\}$. For each $\Gamma^k_2$ consider the ideal $I'_k$ obtained from the intersection of $I$ with ${\bb C}[z_i\colon \ i\in \Gamma^k_2]$. Again, $I$ is generated by ${\bb C}[z_i\colon\ i\in\Gamma_1]$ and the collection $\bigcup\limits_k I'_k$ and thus the zero variety $Z(I)$ is the intersection of the zero variety of ${\bb C}[z_i\colon \ i\in \Gamma_1]$ and the collection $\{Z(I'_k)\}$. Since $I'_k\ne {\bb C}[z_i\colon \ i\in \Gamma^k_2]$, we see that the dimension of $Z(I'_k)$ is strictly less than the cardinality of $\Gamma^\alpha_2$. Since the dimension of $Z(I)$ is one, we see that there can only be one element in the partition of $\Gamma_2$. Hence we can choose an $i_0$ in $\Gamma_2$ and define continuous, monotonically increasing functions $\{\psi_i\}^m_{i=1}$ from $[0,\infty)$ to $[0,\infty)$ so that for $\pmb{z}$ in $Z(I)$ we have $|z_i| = \psi_i(|z_{i_0}|)$ for $1\le i\le m$. Note that $\psi_i\equiv 0$ for $i$ in $\Gamma_1$ and $\psi_{i_0}(x) = x$.

It appears likely that a necessary condition for an ideal $I$ to have an absolutely determining zero variety $Z(I)$, is for $Z(I)$ to have dimension one.

Now one can extend the argument from Theorem \ref{thm1} to establish the following result. Note that when $I$ is homogeneous and $\Omega_\varphi$ is symmetic in all the variables, a better result follows from \cite{GW}.

\begin{tm}\label{thm2}
If ${\cl M}$ is the closure of a radical, bivariate, quasi-homogeneous ideal $I$ in ${\bb C}[\pmb{z}]$ with  zero variety $Z(I)$ of dimension one, then ${\cl M}$ is essentially reductive.
\end{tm}

\begin{proof}
Observe that if we consider a linear polynomial $p_{\pmb{a}}(\pmb{z}) = a_1z_1 +\cdots+ a_mz_m$ for $\pmb{a}$ in ${\bb C}^m\backslash \{\pmb{0}\}$ satisfying $Z(I)\cap \partial\Omega_p\cap Z(p_{\pmb{a}}) =\emptyset$, then it will follow from the preceding argument and Lemma \ref{lem2.4} that $B_{p_{\pmb{a}}}$ is compact. Since the dimension of ${Z}(I)$ is 1 and that of ${Z}(p_{\pmb{a}})$ is $m-1$, then the intersection $Z(I)\cap Z(p_{\pmb{a}})$ consists of just the origin $(0,\ldots, 0)$ for a dense open set $U$ of $\pmb{a}$ in ${\bb C}^m$. This implies that $B_{p_{\pmb{a}}}$ is compact for $\pmb{a}$ in $U$ from which the result follows.
\end{proof}

Can we eliminate the hypothesis that $I$ is radical? Just as in the case of ${\bb C}[z_1,z_2]$, there are two steps:\ (1)~~handle the case in which the $I_i$ are generated by a power of $z_i$ and the $I_{i,j}$ by a power of the prime polynomial $(z^s_i-\alpha z^t_j)^v$ and (2)~~handle the case in which one allows $I_{i,j}$ to be generated by a product of powers of prime polynomials $\prod\limits_k (z^{s_k}_i-\alpha_k z^{t_k}_j)^{v_k}$.

One approach to  the first step would be to relate the essential reductivity of a quasi-homogeneous bivariate ideal $I$ with one dimensional zero variety $Z(I)$ to the essential reductivity of its radical $\sqrt I$. But we have made no progress in doing that even under the assumption that each $I_{i,j}$ is generated by a single prime polynomial. For the second step, one could use Lemma \ref{lem2.3a} to reduce to case 1 if one knew under what circumstances the ideal $I\cap J$ is generated by the collection $\{I_{i,j}\cap J_{i,j}\}$.

It is possible that the proof in \cite{GW2} showing that  a quasi-homogeneous submodule in $H^2_2$ is $p$-essentially reductive for $p>2$, also carries over in the generality of Theorem \ref{thm1} but the author has not verified that. It does seem likely that the result on identifying the $K$-homology class in $K_1(\partial \Omega_\varphi\cap Z({\cl M}))$ carries over. One may need an extension of the index theorem for Toeplitz operators on strongly pseudo-convex domains due to Boutet de Monvel \cite{M}  to this more general class of Reinhardt domains.

Finally, it seems likely that this $K$-homology class agrees with the fundamental class defined by $\partial \Omega_\varphi\cap Z({\cl M})$ as conjectured in \cite{Dou3}.

\end{document}